\documentclass[a4paper,11pt]{amsart}
\usepackage[plainpages=false]{hyperref}
\usepackage{amsfonts,latexsym,rawfonts,amsmath,amssymb,amsthm, mathrsfs, lscape}
\usepackage{verbatim}

\usepackage[all]{xy}
\usepackage{graphicx,psfrag}

\usepackage{array, tabularx}

\usepackage{setspace}

\newtheorem{thm}{Theorem}[section]

\newtheorem{prob}[thm]{Question}

\theoremstyle{remark}
\newtheorem{rmk}[thm]{Remark}

\theoremstyle{definition}

\numberwithin{equation}{section}

\def\p{\partial}
\def\R{\mathbb{R}}

\def\i{\sqrt{-1}}
	
\def\D{\Delta}

\def\cF{{\mathcal F}}

\def\cH{{\mathcal H}}

\begin{document}

\title[On the space of Kahler potentials]{On the space of  Kahler potentials}

\author{Weiyong He}

\address{Department of Mathematics, University of Oregon, Eugene, Oregon, 97403}
\email{whe@uoregon.edu}

\begin{abstract}We consider the geodesic equation for the generalized Kahler potential with only mixed second derivatives bounded. We show that given such two generalized Kahler potentials, there is a unique geodesic segment such that for each point on the geodesic, the generalized Kahler potential has uniformly bounded mixed second derivatives (in manifold directions). This generalizes a fundamental theorem of Chen \cite{Chen00} on the space of Kahler potentials. 
\end{abstract}

\maketitle

\section{Introduction}
Let $(M, [\omega])$ be a compact Kahler manifold of complex dimension $n$. The space of Kahler potentials (in the class of $[\omega]$) is given by
\[
\cH=\{\phi\in C^\infty: \omega_\phi=\omega+\i \p\bar\p \phi>0\}. 
\]
Mabuchi \cite{Mabuchi86} defined a natural metric on $\cH$ by, 
for $\psi_1, \psi_2\in T_\phi\cH$,
\[
\langle \psi_1, \psi_2\rangle_\phi=\int_M \psi_1\psi_2 \omega_\phi^n. 
\]
For any path $\phi(t)\in \cH$, then the geodesic equation is given by
\[
\phi_{tt}-|\nabla\phi_t|^2_{\omega_\phi}=0. 
\]
For any interval $I$ in $\R$, denote $U=I\times S^1$. We use $(z, w)$ to denote points on $M\times U$. 
The geodesic equation is equivalent to the homogeneous complex Monge-Ampere equation (assuming for each $w$, $\phi$ defines a strictly positive Kahler metric, see \cite{Semmes}, \cite{Donaldson99}) 
\begin{equation}\label{E-hcma}
\Omega_\phi^{n+1}=(\phi_{tt}-|\nabla\phi_t|^2_\phi)\frac{\omega_\phi^n}{\omega^n}=0,
\end{equation}
where $\Omega_\phi=\pi^{*}\omega_0+\p\bar \p_{w, z} \phi$, $\pi: M\times U\rightarrow M$ is the projection onto $M$ and $\phi$ is regarded as a $S^1$ invariant function on $M\times U$.

In \cite{Donaldson99}, Donaldson proposed a program which tight up the problems in Kahler geometry regarding the canonical metrics with the geometric structure structure of $\cH$, in particular the geodesic equation plays an important role. 
A fundamental result of Chen \cite{Chen00} asserts that for $I=[0, 1]$ and $\phi_0, \phi_1\in \cH$, there exists a unique $C^{1, 1}$ solution of \eqref{E-hcma} in the sense that
\[
\|\phi\|_{C^1}+\max\{|\p\bar\p_{w, z} \phi|\}\leq C.
\]
Chen's result partially answered Donaldson's conjecture \cite{Donaldson99} and proved that $\cH$ is a metric space with Mabuchi's metric. It turns out to be important to consider generalized Kahler potentials in various problems in Kahler geometry. Fixing a background metric $(M, \omega)$, we are mainly interested in the space of 
the generalized $C^{1, 1}$ Kahler potentials as
\[
\cH_{1, 1}=\{\phi: \omega_\phi\geq 0, \|\phi\|_{C^1}<\infty, 0\leq n+\D\phi<\infty\},
\]
and the space of bounded Kahler potentials
\[
\cH_\infty=\{\phi: \omega_\phi\geq 0, \|\phi\|_{L^\infty}<\infty\},
\]
where the positivity of $\omega_\phi=\omega+\i \p\bar \p \phi$ is understood in the sense of currents. Recently Berndtsson \cite{Bern11} observed that for $\phi_0, \phi_1\in \cH_\infty$, there is a unique generalized solution of \eqref{E-hcma} which is uniformly bounded and moreover $\phi_t$ is also uniformly bounded. Hence in this case the distance function in $\cH_\infty$ is actually well-defined but it remains unclear that whether $\cH_\infty$ is a metric space. The point is that it is not clear that the geodesic minimizes the distance among all curves (which is equivalent to the triangle inequality). In \cite{Chen00} the nondegeneracy of Kahler metrics on two end points $\phi_0, \phi_1$ is crucial for triangle inequality; it is also important to notice that the boundary estimate of $\phi_{tt}$ also depends on the nondegeneracy of two end points (as well as smoothness). It is actually an interesting question to understand the structure of $\cH_{1, 1}$ (and/or $\cH_\infty$) and its relation with the metric completion of $\cH$, for example. 

The main result is concerned with the regularity of the geodesic segment when $\phi_0, \phi_1\in \cH_{1, 1}$.  

\begin{thm}\label{T-main}For any $\phi_0, \phi_1\in \cH_{1, 1}$, there exists a unique generalized solution of \eqref{E-hcma} such that for each $t\in [0, 1]$, 
\[
0\leq n+\D\phi(t)\leq C,
\]
where $C$ is a uniform constant depending only on $n, \|\phi_i\|_{L^\infty}$ and $\sup \D\phi_i$, $i=1, 2$ and the geometry of $(M, \omega)$.
\end{thm}

\begin{rmk}Theorem \ref{T-main} answers a problem of X.-X. Chen \cite{Chen}. 
\end{rmk}

To prove Theorem \ref{T-main}, the new ingredient is the following a priori estimates,

\begin{thm}\label{T-estimate}Let $\phi_0, \phi_1\in \cH$, and let  $f$ be a smooth function and let $\epsilon$ be a positive constant. There is a unique smooth solution of the equation
\begin{equation}\label{E-hcma1}
(\phi_{tt}-|\nabla\phi_t|^2_\phi)\frac{\omega_\phi^n}{\omega^n}=\epsilon e^f, \phi(0, \cdot)=\phi_0, \phi(1, \cdot)=\phi_1.
\end{equation}
Moreover there is a uniform constant $C$ depending only on $n, \|\phi_i\|_{L^\infty}, \sup\D\phi_i,$ $\sup f,$ $\sup \epsilon$, $\|\nabla f\|_{L^\infty}$, 
$\inf \D f$ and the geometry of $(M, g)$ (the upper bound of scalar curvature and the lower bound of the bisectional curvature) such that for each $t\in [0, 1]$,
\[
\|\phi\|_{C^1}+\D\phi\leq C. 
\]
\end{thm}

When $\phi_0, \phi_1\in \cH$, Chen \cite{Chen00} proved the existence of  the unique smooth solution of \eqref{E-hcma1}  and moreover,  there exists a uniform constant $\tilde C$, 
\[
\|\phi\|_{C^1}+\phi_{tt}+\D\phi\leq \tilde C,
\]
where $\tilde C$ does not depend on the lower bound of $\epsilon$. Chen \cite{Chen00} obtained the generalized solution of \eqref{E-hcma} with the the above uniform estimate by letting $\epsilon\rightarrow 0$;  however it is important to note that in his result, $\tilde C$ depends on the strictly positive lower bound of $n+\D\phi_0$ and $n+\D\phi_1$ in a crucial way. Hence the only new ingredient in Theorem \ref{T-estimate} is the estimate of $\D\phi$, independent of the strictly lower bound of $n+\D\phi_0, n+\D\phi_1$.

To establish such an estimate in Theorem \ref{T-estimate},   we follow the similar strategy of Yau's $C^2$ estimate \cite{Yau78} of complex Monge-Ampere equation but we will deal with \eqref{E-hcma1} directly. Our treatment here takes advantage of the product structure of $M\times [0, 1]$ and hence one can separate the manifold direction and $t$ direction respectively. This actually makes possible to estimate directly $\D\phi$ instead of $\phi_{tt}+\D\phi$. A similar situation is already considered in \cite{ChenHe} (see Lemma 3.9, the Ricci curvature assumption there is only technical and not essential) to deal with Donaldson's equation \cite{Donaldson07}. We also note that our computation works for general right hand side with slight modification; for simplicity, we will choose the right hand side to be the special case of $\epsilon e^f$ and it is sufficient for geometric applications. 

We should mention that the study of complex Monge-Ampere equation (nondegenerate/degenerate) has a long history in Kahler geometry and in complex analysis. It has been studied extensively and it is still very active in various settings. Hence we will refer readers to, for example, \cite{K2000} for more references.\\

\noindent{\bf Acknowledgement: } The problem considered here is motivated in part by a question in a recent preprint \cite{He2012} and the author would like to thank Song Sun  for numerous enlightening discussions  and suggestions during the preparation of these two papers. The author is  grateful to Professor X.-X. Chen for constant support and encouragements; our joint work \cite{ChenHe} has definite influence on the current work.

\section{The a priori estimate}

In this section we prove Theorem \ref{T-estimate}. Theorem \ref{T-main} is a rather straightforward consequence of Theorem \ref{T-estimate}, where in particular the uniqueness of generalized solution in the sense of Bedford-Taylor \cite{bt01} will be used. For Theorem \ref{T-estimate}, to solve the equation one needs to establish the a priori estimates and use continuity method, and
we refer to \cite{Chen00} for the details. In particular, the $C^0$ estimate is a standard maximum principle argument. The gradient estimate was first proved in \cite{Chen00} by a blowing-up argument when $\phi_0, \phi_1\in \cH$ and one considers the $C^2$ estimates first. A more direct gradient estimate is established in \cite{Blocki09}. Hence we will only consider the only new ingredient -- the a priori estimate of $\D\phi$, independent of strictly positive bound of $n+\D\phi_i$, $i=1, 2$.

\begin{proof}
We consider the equation 
\begin{equation}
(\phi_{tt}-|\nabla \phi_t|^2_\phi) \det(g_{i\bar j}+\phi_{i\bar j})=\epsilon e^f \det(g_{i\bar j}),
\end{equation}
where we use the complex notation as follows, 
\[
|\nabla \phi|^2=g^{k\bar l}\phi_k\phi_{\bar l}=\phi_k\phi_{\bar k}, |\nabla \phi_t|^2_\phi=g^{i\bar j}_\phi \phi_{ti}\phi_{t\bar j}, \D h=g^{k\bar l}h_{k\bar l}, \D_\phi h=g^{k\bar l}_\phi h_{k\bar l}. 
\]
We rewrite the equation as
\begin{equation}\label{cma-1}
\log(\phi_{tt}-|\nabla \phi_t|^2_\phi)+\log \det(g_{i\bar j}+\phi_{i\bar j})=\log \epsilon+f+\log \det(g_{i\bar j}). 
\end{equation}
The linearization of the left hand side of \eqref{cma-1} is given by
\begin{equation}
\begin{split}
Dh=&\D_\phi h+\frac{h_{tt}+g^{i\bar l}_\phi g^{k\bar j}_\phi h_{k\bar l}\phi_{ti}\phi_{t\bar j}}{\phi_{tt}-|\nabla \phi_t|^2_\phi}\\
&-\frac{g^{i\bar j}_\phi (h_{ti}\phi_{t\bar j}+h_{t\bar j}\phi_{ti})}{\phi_{tt}-|\nabla \phi_t|^2_\phi}. 
\end{split}
\end{equation}
We compute
\begin{equation}\label{E-1-4}
\begin{split}
D(\D \phi)=&\D_\phi (\D \phi)+\frac{\D\phi_{tt}+g^{i\bar l}_\phi g^{k\bar j}_\phi (\D \phi)_{k\bar l}\phi_{ti}\phi_{t\bar j}}{\phi_{tt}-|\nabla\phi_t|^2_\phi}\\
&-\frac{g^{i\bar j}_\phi ((\D\phi)_{ti}\phi_{t\bar j}+(\D\phi)_{t\bar j}\phi_{ti})}{\phi_{tt}-|\nabla \phi_t|^2_\phi}. 
\end{split}
\end{equation}
A direct computation by differentiating \eqref{cma-1}, similar as in \cite{Yau78}, gives 
\begin{equation}\label{E-1-5}
\begin{split}
\D_\phi(\D \phi)+\frac{\D\left(\phi_{tt}-|\nabla \phi_t|^2_\phi\right)}{\phi_{tt}-|\nabla\phi_t|^2_\phi}=&g^{i\bar j}_\phi g^{p\bar q}_\phi \phi_{i\bar q k}\phi_{p\bar j\bar k}+\frac{\left|\nabla (\phi_{tt}-|\nabla\phi_t|^2_\phi)\right|^2}{\left(\phi_{tt}-|\nabla\phi_t|^2_\phi\right)^2}\\
&+\D f+I,
\end{split}
\end{equation}
where 
\[
I=\sum_{i, j, k, l}\left(g^{i\bar j}_\phi R_{i\bar j l\bar l}-R_{i\bar i l\bar l}+g^{i\bar j}_\phi R_{i\bar j k\bar l}\phi_{l\bar k}\right). 
\]

Let the bisectional curvature of $(M, g)$ satisfie, for some nonnegative number $B$
\[
R_{i\bar jk\bar l}\geq -B(g_{i\bar j}g_{k\bar l}+g_{i\bar l}g_{k\bar j})
\]
One can estimate, for example see \cite{Yau78}, 
\begin{equation}\label{E-yau}
\begin{split}
I=&g^{i\bar j}_\phi R_{i\bar j k\bar l} (g_{l\bar k}+\phi_{l\bar k})-R\\
=&\sum_{i, k}\frac{1+\phi_{k\bar k}}{1+\phi_{i\bar i}} R_{i\bar i k\bar k}-R\\
\geq &-B(g_{i\bar i}g_{k\bar k}+g_{i\bar k}g_{k\bar i})\frac{1+\phi_{k\bar k}}{1+\phi_{i\bar i}}-R\\
= &-B\sum_{i, k}\frac{1+\phi_{k\bar k}}{1+\phi_{i\bar i}}-B-R\\
=&-\sum_i\frac{B(n+\D\phi)}{1+\phi_{i\bar i}}-B-R.
\end{split}
\end{equation}
We compute, by assuming $g_{i\bar j}=\delta_{ij}, \p g_{i\bar j}=0$ at the point $P$, 
\begin{equation}\label{E-1-6}
\begin{split}
\D(|\nabla \phi_{t}|_\phi^2)=&g^{k\bar l}\frac{\p^2}{\p z_k\p z_{\bar l}}\left(g^{i\bar j}_\phi \phi_{ti}\phi_{t\bar j}\right)\\
=&g^{k\bar l}\p_k\left(\frac{\p g^{i\bar j}_\phi}{\p z_{\bar l}}\phi_{ti}\phi_{t\bar j}+g^{i\bar j}_\phi\p_{\bar l}\phi_{ti} \phi_{t\bar j}+g^{i\bar j}_\phi\phi_{ti}\p_{\bar l}\phi_{t\bar j}\right)\\
=&g^{k\bar l}\left(\frac{\p^2 g^{i\bar j}_\phi}{\p z_k \p z_{\bar l}} \phi_{ti}\phi_{t\bar j}+{\p_{\bar l}g^{i\bar j}_\phi} \p_k\phi_{ti}\phi_{t\bar j}+{\p_{\bar l}g^{i\bar j}_\phi} \phi_{ti} \phi_{tk\bar j}\right)\\
&+g^{k\bar l}\left(\p_k g^{i\bar j}_\phi\phi_{ti\bar l}\phi_{t\bar j}+g^{i\bar j}_\phi \p_k\phi_{ti\bar l}\phi_{t\bar j}+g^{i\bar j}_\phi\phi_{ti\bar l}\phi_{tk\bar j}\right)\\
&+g^{k\bar l}\left(\p_k g^{i\bar j}_\phi\phi_{ti}\p_{\bar l}\phi_{t\bar j}+g^{i\bar j}_\phi\p_k\phi_{ti}\p_{\bar l}\phi_{t\bar j}+g^{i\bar j}_\phi \phi_{ti}\p^2_{k, \bar l}\phi_{t\bar j}\right)
\end{split}
\end{equation}
At the point $P$, we compute
\begin{equation}\label{E-1-7}
\begin{split}
\p_k g^{i\bar j}_\phi=-&g^{i\bar q}_\phi g^{p\bar j}_\phi \p_k\left(g_{p\bar q}+\phi_{p\bar q}\right)=-g^{i\bar q}_\phi g^{p\bar j}_\phi\phi_{p\bar q k}\\
\p_{\bar l} g^{i\bar j}_\phi=-&g^{i\bar q}_\phi g^{p\bar j}_\phi \p_{\bar l}\left(g_{p\bar q}+\phi_{p\bar q}\right)=-g^{i\bar q}_\phi g^{p\bar j}_\phi\phi_{p\bar q \bar l}
\end{split}
\end{equation}
and hence we have
\begin{equation}\label{E-1-8}
\begin{split}
g^{k\bar l}\p^2_{k, \bar l}\left(g^{i\bar j}_\phi\right)=&-g^{k\bar l} g^{i\bar q}_\phi g^{p\bar j}_\phi \p^{2}_{k, \bar l}(g_{p\bar q}+\phi_{p\bar q})\\
&+g^{k\bar l}g^{i\bar s}_\phi g^{r\bar q}_\phi g^{p\bar j}_\phi\p_k(g_{r\bar s}+\phi_{r\bar s})\p_{\bar l}(g_{p\bar q}+\phi_{p\bar q})\\
&+g^{k\bar l}g^{i\bar q}_\phi g^{p\bar s}_\phi g^{r\bar j}_\phi \p_k(g_{r\bar s}+\phi_{r\bar s})\p_{\bar l}(g_{p\bar q}+\phi_{p\bar q})\\
=&-g^{k\bar l} g^{i\bar q}_\phi g^{p\bar j}_\phi \p^{2}_{k, \bar l}(g_{p\bar q})\\
&-g^{i\bar q}_\phi g^{p\bar j}_\phi \p^{2}_{p, \bar q}\left(g^{k\bar l}\p^2_{k, \bar l}\phi\right)+g^{i\bar q}_\phi g^{p\bar j}_\phi \p^2_{p, \bar q}\left(g^{k\bar l}\right)\phi_{k\bar l}\\
&+g^{k\bar l}g^{i\bar s}_\phi g^{r\bar q}_\phi g^{p\bar j}_\phi\p_k(g_{r\bar s}+\phi_{r\bar s})\p_{\bar l}(g_{p\bar q}+\phi_{p\bar q})\\
&+g^{k\bar l}g^{i\bar q}_\phi g^{p\bar s}_\phi g^{r\bar j}_\phi \p_k(g_{r\bar s}+\phi_{r\bar s})\p_{\bar l}(g_{p\bar q}+\phi_{p\bar q})\\
=&g^{k\bar l}g^{i\bar q}_\phi g^{p\bar j}_\phi R_{p\bar q k\bar l}-g^{i\bar q}_\phi g^{p\bar j}_\phi (\D\phi)_{p\bar q}+g^{i\bar q}_\phi g^{p\bar j}_\phi R_{k\bar l p\bar q} \phi_{k\bar l}\\
&+g^{k\bar l}g^{i\bar s}_\phi g^{r\bar q}_\phi g^{p\bar j}_\phi\left(\phi_{r\bar sk}\phi_{p\bar q \bar l}+\phi_{p\bar qk}\phi_{r\bar s\bar l}\right)
\end{split}
\end{equation}
By \eqref{E-1-6}, \eqref{E-1-7} and\eqref{E-1-8}, we compute,
\begin{equation}\label{E-1-9}
\begin{split}
\D(|\nabla \phi_{t}|_\phi^2)=&g^{k\bar l}g^{i\bar q}_\phi g^{p\bar j}_\phi R_{p\bar q k\bar l}\phi_{ti}\phi_{t\bar j}-g^{i\bar q}_\phi g^{p\bar j}_\phi (\D\phi)_{p\bar q}\phi_{ti}\phi_{t\bar j}\\
&+g^{i\bar q}_\phi g^{p\bar j}_\phi R_{k\bar l p\bar q} \phi_{l\bar k}\phi_{ti}\phi_{t\bar j}\\
&+g^{k\bar l}g^{i\bar s}_\phi g^{r\bar q}_\phi g^{p\bar j}_\phi\left(\phi_{r\bar sk}\phi_{p\bar q \bar l}+\phi_{p\bar qk}\phi_{r\bar s\bar l}\right)\phi_{ti}\phi_{t\bar j}\\
&-g^{k\bar l}g^{i\bar q}_\phi g^{p\bar j}_\phi \phi_{p\bar q\bar l}(\phi_{tik}\phi_{t\bar j}+\phi_{ti}\phi_{tk\bar j})\\
&-g^{k\bar l}g^{i\bar q}_\phi g^{p\bar j}_\phi\phi_{p\bar q k}(\phi_{ti\bar l}\phi_{t\bar j}+\phi_{ti}\phi_{t\bar j\bar l})\\
&+g^{i\bar j}_\phi (\D\phi)_{ti}\phi_{t\bar j}+g^{i\bar j}_\phi\phi_{ti}(\D \phi)_{t\bar j}\\
&+g^{k\bar l}g^{i\bar j}_\phi (\phi_{ti\bar l}\phi_{t\bar j k}+\phi_{tik}\phi_{t\bar j\bar l})
\end{split}
\end{equation}
Hence we compute, by \eqref{E-1-4}, \eqref{E-1-5} and \eqref{E-1-9},
\begin{equation}\label{E-1-11}
\begin{split}
D(\D \phi)=&\D f+g^{i\bar j}_\phi g^{p\bar q}_\phi \phi_{i\bar q k}\phi_{p\bar j\bar k}+\frac{\left|\nabla (\phi_{tt}-|\nabla\phi_t|^2_\phi)\right|^2}{\left(\phi_{tt}-|\nabla\phi_t|^2_\phi\right)^2}\\
&+I+\frac{II+III}{\phi_{tt}-|\nabla \phi_t|^2_\phi}
\end{split}
\end{equation}
where
\[
II=g^{k\bar l}g^{i\bar q}_\phi g^{p\bar j}_\phi R_{p\bar q k\bar l}\phi_{ti}\phi_{t\bar j}+g^{i\bar q}_\phi g^{p\bar j}_\phi R_{k\bar l p\bar q} \phi_{l\bar k}\phi_{ti}\phi_{t\bar j}
\]
and 
\begin{equation}\label{E-III}
\begin{split}
III=&g^{k\bar l}g^{i\bar s}_\phi g^{r\bar q}_\phi g^{p\bar j}_\phi\left(\phi_{r\bar sk}\phi_{p\bar q \bar l}+\phi_{p\bar qk}\phi_{r\bar s\bar l}\right)\phi_{ti}\phi_{t\bar j}\\
&-g^{k\bar l}g^{i\bar q}_\phi g^{p\bar j}_\phi \phi_{p\bar q\bar l}(\phi_{tik}\phi_{t\bar j}+\phi_{ti}\phi_{tk\bar j})\\
&-g^{k\bar l}g^{i\bar q}_\phi g^{p\bar j}_\phi\phi_{p\bar q k}(\phi_{ti\bar l}\phi_{t\bar j}+\phi_{ti}\phi_{t\bar j\bar l})\\
&+g^{k\bar l}g^{i\bar j}_\phi (\phi_{ti\bar l}\phi_{t\bar j k}+\phi_{tik}\phi_{t\bar j\bar l}).
\end{split}
\end{equation}

Note that as an Hermitian matrix, $(\phi_{ti}\phi_{t\bar j})$ is nonnegative definite and hence
we can then estimate, where we choose a coordinate at the point such that $g_{i\bar j}=\delta_{ij}, \phi_{i\bar j}=\phi_{i\bar i}\delta_{ij}$, 
\begin{equation}\label{E-1-12}
\begin{split}
II=&g^{i\bar q}_\phi g^{p\bar j}_\phi R_{p\bar qk \bar l}(g_{l\bar k}+\phi_{l\bar k})\phi_{ti}\phi_{t\bar j}\\
\geq &-B(g_{p\bar q}g_{k\bar l}+g_{p\bar l}g_{k\bar q}) g^{i\bar q}_\phi g^{p\bar j}_\phi (g_{l\bar k}+\phi_{l\bar k})\phi_{ti}\phi_{t\bar j}\\
=&-B(n+\D \phi)\sum_{i}\frac{\phi_{ti}\phi_{t\bar i}}{(1+\phi_{i\bar i})^2}-B|\nabla \phi_t|_\phi^2\\
>& -2B(n+\D \phi)\sum_{i}\frac{\phi_{ti}\phi_{t\bar i}}{(1+\phi_{i\bar i})^2}.
\end{split}
\end{equation}
where we use the observation  for each $i$, 
\[
n+\D\phi> 1+\phi_{i\bar i}, (n+\D\phi)^{-1}<(1+\phi_{i\bar i})^{-1}. 
\]
and hence we can estimate
\[
\frac{|\nabla\phi_t|^2_\phi}{n+\D\phi}=\sum_i\frac{\phi_{ti}\phi_{t\bar i}}{1+\phi_{i\bar i}}(n+\D\phi)^{-1}<\sum_{i}\frac{\phi_{ti}\phi_{t\bar i}}{(1+\phi_{i\bar i})^2}.
\]
It then follows from \eqref{E-1-11}, using \eqref{E-yau} and \eqref{E-1-12} 
\begin{equation}\label{E-1-14}
\begin{split}
D(\D\phi)\geq &\D f+g^{i\bar j}_\phi g^{p\bar q}_\phi \phi_{i\bar q k}\phi_{p\bar j\bar k}+\frac{III}{\phi_{tt}-|\nabla\phi_t|^2_\phi}\\
&-B\sum_{i}\frac{n+\D\phi}{1+\phi_{i\bar i}}-B-R\\
&-2B\frac{(n+\D\phi)}{\phi_{tt}-|\nabla\phi_t|^2_\phi} \sum_i\frac{\phi_{ti}\phi_{t\bar i}}{(1+\phi_{i\bar i})^2}
\end{split}
\end{equation}
A direct computation gives, for $h>0$,
\begin{equation}\label{E-1-15}
\begin{split}
D(\log h)=&\frac{D h}{h}-\frac{g^{k\bar l}_\phi h_kh_{\bar l}}{h^2}-\frac{\left(h_t-g^{i\bar l}_\phi\phi_{ti}h_{\bar l}\right)\left(h_t-g^{k\bar j}_\phi \phi_{t\bar j}h_k\right)}{h^2(\phi_{tt}-|\nabla\phi_t|^2_\phi)}
\end{split}
\end{equation}
We also have
\begin{equation}\label{E-1-16}
\begin{split}
D\phi&=(n+1)-g^{i\bar j}_\phi g_{i\bar j}-\frac{g^{i\bar l}_\phi g^{k\bar j}_\phi g_{k\bar l}\phi_{ti}\phi_{t\bar j}}{\phi_{tt}-|\nabla \phi_t|^2_\phi}\\
&=(n+1)-\sum_i\frac{1}{1+\phi_{i\bar i}}-\frac{1}{\phi_{tt}-|\nabla\phi_t|^2_\phi} \sum_i\frac{\phi_{ti}\phi_{t\bar i}}{(1+\phi_{i\bar i})^2}
\end{split}
\end{equation}
Take $h=n+\D\phi$, and let
\[
A=\frac{\left(h_t-g^{i\bar l}_\phi\phi_{ti}h_{\bar l}\right)\left(h_t-g^{k\bar j}_\phi \phi_{t\bar j}h_k\right)}{h^2(\phi_{tt}-|\nabla\phi_t|^2_\phi)}
\]
Then we compute, using \eqref{E-1-15} and \eqref{E-1-16}
\begin{equation}\label{E-1-17}
\begin{split}
D(\log h-C\phi)=&\frac{D(\D\phi)}{n+\D\phi}-\frac{g^{k\bar l}_\phi h_kh_{\bar l}}{h^2}-A-(n+1)C\\
&+\sum_i\frac{C}{1+\phi_{i\bar i}}+\frac{C}{\phi_{tt}-|\nabla\phi_t|^2_\phi} \sum_i\frac{\phi_{ti}\phi_{t\bar i}}{(1+\phi_{i\bar i})^2}
\end{split}
\end{equation}
Applying \eqref{E-1-14}, it follows from \eqref{E-1-17} that
\begin{equation}\label{E-1-18}
\begin{split}
D(\log h-C\phi)\geq &\frac{\D f-B-R}{n+\D\phi}+\frac{g^{i\bar j}_\phi g^{p\bar q}_\phi \phi_{i\bar q k}\phi_{p\bar j\bar k}}{n+\D\phi}-\frac{g^{k\bar l}_\phi h_kh_{\bar l}}{h^2}\\
&+\sum_i\frac{C-B}{1+\phi_{i\bar i}}-(n+1)C\\
&+\frac{(C-2B) \sum_i\frac{\phi_{ti}\phi_{t\bar i}}{(1+\phi_{i\bar i})^2}}{\phi_{tt}-|\nabla\phi_t|^2_\phi}\\
&+\frac{III(n+\D\phi)^{-1}}{\phi_{tt}-|\nabla\phi_t|^2_\phi}-A
\end{split}
\end{equation}
We have the estimate (see \cite{Yau78} (2.15) for example),
\[
g^{i\bar j}_\phi g^{p\bar q}_\phi \phi_{i\bar q k}\phi_{p\bar j\bar k}\geq (n+\D\phi)^{-1} g^{k\bar l}_\phi (\D\phi)_k(\D\phi)_{\bar l}. 
\]
We also have the estimate
\begin{equation}\label{E-1-19}
III\geq A (\phi_{tt}-|\nabla\phi_t|^2_\phi).
\end{equation}
It is a straightforward computation to establish \eqref{E-1-19} and we will derive this in \eqref{E-III1}.
Hence if we choose $C$ sufficiently large (but finite number) such that
\[
C-2B-R+\inf \D f\geq 1, C-2B\geq 1,
\]
we have the estimate
\[
D(\log h-C\phi)\geq \sum_i\frac{1}{1+\phi_{i\bar i}}-(n+1)C.
\]
Observe that
\[
D(t^2)=2(\phi_{tt}-|\nabla \phi_t|^2_\phi)^{-1}.
\]
Hence we have
\[
D(\log h-C\phi+t^2)> \sum_i(1+\phi_{i\bar i})^{-1}+(\phi_{tt}-|\nabla\phi_t|^2_\phi)^{-1}-(n+1)C.
\]
Now following \cite{Yau78} (2.19, 2.20) for example, we can get that
\[
 \sum_i(1+\phi_{i\bar i})^{-1}+(\phi_{tt}-|\nabla\phi_t|^2_\phi)^{-1}\geq (n+\D\phi+\phi_{tt}-|\nabla\phi_t|^2_\phi)^{1/n} \epsilon^{-1/n} e^{\frac{-f}{n}}
\]
Either $\log h-C\phi+t^2$ obtain its maximum interior at some point $P$, then 
at $P$, 
\[
D(\log h-C\phi+t^2)\leq 0,
\]
it then follows that
\[
(n+\D\phi)(P)< C \exp(f/n) \epsilon^{1/n};
\]
for some uniformly bounded constant $C=C(B,  n, \inf \D f)$. Clearly in this case
\begin{equation}
n+\D\phi< C_0 \exp(f/n)\epsilon^{1/n},
\end{equation}
where $C_0=C_0(\|\phi\|_{L^\infty}, B, n, \inf \D f)$. 
Or $\log h-C\phi+t^2$ obtains its maximum on boundary,  then we also have
\[
n+\D\phi\leq C, 
\]
where $C$ depends on $\|\phi\|_{L^\infty}, \sup \D\phi_0$ and $\sup \D\phi_1$. 

\end{proof}

Now we establish \eqref{E-1-19}. We need to show that, 

\begin{equation}\label{E-III1}
III (n+\D\phi)\geq \left((\D\phi)_t-g^{i\bar l}_\phi\phi_{ti}(\D\phi)_{\bar l}\right)\left((\D\phi)_t-g^{k\bar j}_\phi \phi_{t\bar j}(\D\phi)_k\right),
\end{equation}
where $III$ is given by \eqref{E-III}. 
We rewrite
\[
III=A_1+A_2,
\]
where we set
\[
\begin{split}
A_1=&g^{k\bar l}g^{i\bar s}_\phi g^{r\bar q}_\phi g^{p\bar j}_\phi\phi_{p\bar qk}\phi_{r\bar s\bar l}\phi_{ti}\phi_{t\bar j}+g^{k\bar l}g^{i\bar j}_\phi \phi_{ti\bar l}\phi_{t\bar j k}\\
&-g^{k\bar l}g^{i\bar q}_\phi g^{p\bar j}_\phi \phi_{p\bar q\bar l}\phi_{ti}\phi_{tk\bar j}-g^{k\bar l}g^{i\bar q}_\phi g^{p\bar j}_\phi\phi_{p\bar q k}\phi_{ti\bar l}\phi_{t\bar j}\\
\end{split}
\]
and
\[
\begin{split}
A_2=&g^{k\bar l}g^{i\bar s}_\phi g^{r\bar q}_\phi g^{p\bar j}_\phi\phi_{r\bar sk}\phi_{p\bar q\bar l}\phi_{ti}\phi_{t\bar j}+g^{k\bar l}g^{i\bar j}_\phi\phi_{tik}\phi_{t\bar j\bar l}\\
&-g^{k\bar l}g^{i\bar q}_\phi g^{p\bar j}_\phi \phi_{p\bar q\bar l}\phi_{tik}\phi_{t\bar j}-g^{k\bar l}g^{i\bar q}_\phi g^{p\bar j}_\phi\phi_{p\bar q k}\phi_{ti}\phi_{t\bar j\bar l}
\end{split}
\]
We can estimate
\[
\begin{split}
A_1=&\sum_{r, k}\frac{1}{1+\phi_{r\bar r}}\left\{\left(\sum_i\frac{\phi_{ti}\phi_{r\bar i \bar k}}{1+\phi_{i\bar i}}\right)\left(\sum_j\frac{\phi_{t\bar j}\phi_{j\bar r\bar k}}{1+\phi_{r\bar r}}\right)+\phi_{tr\bar k}\phi_{t\bar r k}\right\}\\
&-\sum_{r, k}\frac{1}{1+\phi_{r\bar r}}\left(\phi_{t\bar r k}\sum_i\frac{\phi_{ti}\phi_{r\bar i\bar k}}{1+\phi_{i\bar i}}+\phi_{tr\bar k}\sum_j\frac{\phi_{t\bar j}\phi_{j\bar r  k}}{1+\phi_{j\bar j}}\right)\\
=&\sum_{r, k}\frac{1}{1+\phi_{r\bar r}} (M_{r k}-\phi_{tr\bar k})(\overline{M_{ r k}}-\phi_{t\bar r k} )\\
\geq & \sum_{r}\frac{1}{1+\phi_{r\bar r}} (M_{r r}-\phi_{tr\bar r})\overline{(M_{r r}-\phi_{tr\bar r})}
\end{split}
\]
where we denote $M_{r k}$ as the matrix (of complex number),
\[
M_{rk}=\sum_i\frac{\phi_{ti}\phi_{r\bar i\bar k}}{1+\phi_{i\bar i}},\; \overline{M_{rk}}=\sum_j\frac{\phi_{t\bar j}\phi_{j\bar r k}}{1+\phi_{j\bar j}}
\]
In particular $A_1\geq 0$ and it follows similarly that $A_2\geq 0$. Moreover, we can estimate
\[
\begin{split}
A_1 (n+\D\phi)\geq& \left(\sum_{r}\frac{1}{1+\phi_{r\bar r}} (M_{r r}-\phi_{tr\bar r})\overline{(M_{rr}-\phi_{tr\bar r})}\right)\left(\sum_r (1+\phi_{r\bar r})\right)\\
\geq & \left(\sum_{r} (M_{rr}-\phi_{tr\bar r})\right)\overline{\left(\sum_r(M_{rr}-\phi_{t\bar r r})\right)}\\
=&\left(\sum_i\frac{\phi_{ti}(\D\phi)_{\bar i}}{1+\phi_{i\bar i}}-\D\phi_t\right)\left(\sum_j\frac{\phi_{t\bar j}(\D\phi)_j}{1+\phi_{j\bar j}}-\D\phi_t\right)\\
=&\left((\D\phi)_t-g^{i\bar l}_\phi\phi_{ti}(\D\phi)_{\bar l}\right)\left((\D\phi)_t-g^{k\bar j}_\phi \phi_{t\bar j}(\D\phi)_k\right).
\end{split}
\]
It completes the proof of \eqref{E-III1}. 

\begin{rmk}If the righthand side is a positive function $F$, then a slight modification can get that
\[
D(log (n+\D\phi)-C\phi+t^2)>\frac{\D(\log F)}{n+\D\phi}+F^{-1/n} (n+\D\phi+\phi_{tt}-|\nabla\phi_t|^2_\phi)-C_1,
\]
where $C_1=C_1(n, B)$. 

\end{rmk}

\section{Discussions}

To prove Theorem \ref{T-main}, we assume $\phi_0, \phi_1\in \cH_{1, 1}$. Then there exists a unique generalized solution $\phi(t)$ of \eqref{E-hcma} with uniformly bounded $\|\phi\|_{L^\infty}$, $|\phi_t|$ and $|\nabla \phi|$. 
We choose a sequence of functions $\phi_i^k\in \cH$, which converges to $\phi_i$ in $C^{1, \alpha}$ and such that $\D\phi_i^k\leq C$. By Theorem \ref{T-estimate}, there exists a sequence of solutions $\phi^k(t)$ of 
\[
(\phi_{tt}-|\nabla\phi_t|^2_\phi)\omega^n_\phi (\omega^n)^{-1}=\frac{1}{k} e^f, \phi^k(0)=\phi_0^k, \phi^k(1)=\phi^k_1. 
\] 
By passing to a subsequence, $\phi^k(t)\rightarrow \phi(t)$ in $C^\alpha$ by the uniqueness of the generalized solution. Moreover by Theorem \ref{T-estimate},  $\D\phi^k\leq C$, where $C$ is independent of the strictly positive lower bound of $n+\D\phi^k_i$, $i=1, 2$. Hence let $k\rightarrow \infty$, we can get that $\D\phi (t)\leq C$.

As mentioned in the introduction, the problem considered here is motivated by the questions in \cite{He2012},

\begin{prob}Is $\cH_{1, 1}$ (or $\cH_\infty$) a metric space?
If the answer is affirmative, does $\cH_{1, 1}$ (or $\cH_\infty$) have nonpositive curvature in the sense of Alexanderov as $\cH$ does (see \cite{CalabiChen})?  What is the relation between $\cH_{1, 1}$ (or $\cH_\infty$) and the metric completion of $\cH$? 
\end{prob}

However, our estimates do not give any information on $\phi_{tt}$ when $\phi_i\in \cH_{1, 1}$; our results seem not to be able to answer the above question directly even though we believe it would be an interesting problem.

\end{document}